\documentclass[a4paper,11pt]{elsarticle}

\usepackage[utf8]{inputenc}

\usepackage{times} 
\usepackage[T1]{fontenc}         
\usepackage{pslatex}
\usepackage{amsmath}
\usepackage{mathtools}
\usepackage{microtype} 
\usepackage{float}
\usepackage[percent]{overpic}
\usepackage{subfig}
\usepackage[format=plain,labelfont=bf]{caption}
\usepackage{xcolor}

\usepackage{tikz}
\usetikzlibrary{shapes.geometric, arrows, shadows,positioning}

\tikzstyle{startstop} = [rectangle, rounded corners, minimum width=3cm, minimum height=1cm,text centered, text width=3cm, draw=black]
\tikzstyle{process} = [rectangle, minimum width=3cm, minimum height=1cm, text centered, text width=4cm, draw=black]
\tikzstyle{decision} = [diamond, minimum width=2cm, minimum height=0.5cm, text centered, draw=black]
\tikzstyle{arrow} = [thick,->,>=stealth]

\title{A novel method for robust and efficient prediction of ice shedding from rotorcraft blades}

\author{Andrea Rausa, Myles Morelli, Alberto Guardone \\[0.5cm]
Politecnico di Milano, \\
Department of Aerospace Science and Technology \\
Via La Masa 34, 20156 Milano, Italy}

\date{}

\begin{document}

\begin{abstract}

\noindent
The safety of rotorcraft operating in cold environments is jeopardised by the possibility of ice accretion on the rotor blades.
Eventually, ice can shed from the blade due to the high centrifugal forces and impact other parts of the rotorcraft or unbalance the rotor. To establish the shedding time and location for rotorcraft, a robust and efficent numerical multi-step icing simulations tool is presented here for predicting the shedding phenomenon. A volume-mesh based approach is used to allow for representing the ice shape. In order to increase the robustness of the method, an interpolation procedure is implemented which establishes the possible occurrence of the shedding event and restricts the search domain. Ice shapes along the blades are computed by means of two-dimensional ice accretion simulations: ice shapes are then interpolated over the blade span. 

\noindent
Numerical results compares fairly well, in terms of shedding time and location, to the experimental ones obtained in the AERTS test facility, thus demonstrating the soundness of the present approach.

\end{abstract}
\maketitle

\thispagestyle{empty}



\section{Introduction}\label{ch:Intro}

In-flight ice accretion poses major concerns regarding the safety of air transportation \cite{cebeci2003aircraft}. Ice accretion on aerodynamic surfaces, like wings and empennage, drastically changes their shape and increases the weight of the aircraft, thus leading to a high degradation of aerodynamic performances.

\vspace{\baselineskip}
\noindent
Ice accretion is not limited to fixed-wing aircraft, but it affects also rotating aerodynamic devices, such as aero-engine compressor or rotorcraft. Ice accretion modifies the lift, drag and pitching moment characteristics of the blade sections. Degradation of the lift characteristics results in the increase of collective input to the main rotor, with consequent increased power requirement for a given flight condition, together with erosion of blade stall margins \cite{gent2000aircraft}. 

\vspace{\baselineskip}
\noindent
Ice shedding from rotor blades is caused by centrifugal forces because of the high rotational velocity. As the ice accretes, centrifugal forces eventually overcome both forces of bonds between ice and blade and within ice itself, causing it to shed. The shedding of ice chunks from the main rotor can be an asymmetric phenomenon, thus leading to strong rotor imbalances and subsequent severe vibrations. Moreover, the shed chunk of ice can represent a ballistic danger since it can hit other parts of the rotorcraft, thus leading to severe limitations to its manoeuvrability. 

\vspace{\baselineskip}

Previous research efforts tried various approach to tackle ice shedding from rotor blades. 

\vspace{\baselineskip}
\noindent
Zhang and Habashi \cite{zhang2012fem} developed a 3D Finite Element Analysis for in-flight ice break-up of rotorcraft. Aerodynamic forces are supposed to be negligible with respect to the centrifugal ones. Shear stress, principal stress as well as Von-Mises stresses are checked at each crack propagation step. If either shear stress or Von-Mises stresses reaches a critical value, the interface bond between ice and airfoil breaks. Meanwhile, if the maximum principal stress has exceeded the cohesive strength of the ice, a crack is initiated. This method is accurate but quite computationally expensive.

\vspace{\baselineskip}
\noindent

Chen and Fu \cite{chen2015numerical} developed an ice shedding model by taking into account the coupling of the failure of the interface between ice and rotor blade surface and the failure of ice itself. A bilinear cohesive zone model (CZM) was introduced to simulate the initiation and propagation of ice/blade interface crack and a maximum stress criterion was used to describe the failure occurred in the ice. This model was tested against experimental results, confirming the accuracy of a FEM structural analysis in describing the ice shedding phenomenon.

\vspace{\baselineskip}
\noindent
Brouwers et al. \cite{brouwers2010experimental} analysed ice shedding from rotor blades. Their approach combines a BEMT (Blade Element Momentum Theory) analysis with LEWICE \cite{wright2005validation} ice accretion code in order to compute ice shapes at selected stations along the blade. Then, starting from the tip and going to the root of the blade, for each station both cohesion and adhesion surfaces are computed, as well as centrifugal forces. The aerodynamic forces were considered negligible for simplicity. Even if such a simplification was made, their method proved to be computationally inexpensive due to the presence of a BEMT model instead of a full 3D flow simulation.

\vspace{\baselineskip}
\noindent
In the previous mentioned works,aerodynamic forces are neglected. This is possible since the tip Mach number of the analysed rotor is usually less than 0.45. This value resulted from a FEM analysis performed by Scavuzzo et al. \cite{scavuzzo1994influence}. In their research, they used a 2D Navier-Stokes solver in order to obtain the flow around an iced airfoil. From this computation, the pressure coefficients was extracted and used to compute the pressure acting on the ice shape. It was found that stresses produced at Mach numbers below 0.45 were not significant enough to contribute to shedding. However they emphasize the fact that in case of high angle of attack or higher Mach numbers, aerodynamic forces cannot be neglected.

\vspace{\baselineskip}
\noindent

In this paper, the ice shedding model used by Brouwers et al. \cite{brouwers2010experimental} is improved by analysing the complete ice shape through the creation of a 3D volume mesh. Through this approach it is possible to partition the ice volume through cutting planes and compute "exactly" the values of centrifugal, adhesion and cohesion forces for each ice piece analysed. Thus it does not require a high discretization of the ice geometry along the span of the blade, resulting in a decrease of the computational cost with respect to the previous works.

\vspace{\baselineskip}
\noindent

The paper is organized as the following. In Section 2 the proposed methodology will be introduced by explaining its main components: the \textit{Iterative cutting} and the \textit{Force-fitting} procedures. In Section 3 the methodolgy will be validated against experimental results. At last, in Section 4 conclusions will be provided as well as new ideas for future developments.
\section{Proposed methodology}

Ice shedding on rotational frames is caused by the centrifugal forces that act on the ice geometry. Shedding occurs if centrifugal forces overcome both the adhesion forces between the ice and the blade and the cohesion forces within the ice. 

The in-house ice accretion tool PoliMIce \cite{gori2015polimice} is used here to simulate the ice accretion. In the accretion process the time scale of the accretion of the ice is separated from the one of the adaptation of the aerodynamic to the new shape. Thus it is possible to operate a multi-step analysis in which during the accretion step the flow properties are frozen.

Through PoliMIce 2D ice shapes are computed at selected stations along the blade. The 3D ice shape is obtained through constant-area extrusion of the aforementioned 2D geometries along the radial direction. The inner volume-mesh is created by the software GMSH \cite{geuzaine2009gmsh} through a Delaunay triangulation algorithm.
The numerical method here developed is divided into two parts: the \textit{Iterative cutting} and the \textit{Force-fitting procedure}.

\subsection{Iterative cutting}\label{ch:Method:Iterative}

\begin{figure}[h!]
	\centering	
	\subfloat[][\label{fig:cutting:First}]{\includegraphics[width=0.49\textwidth]{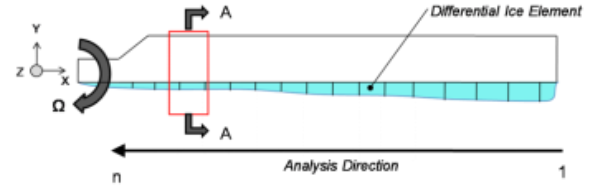}}
	\hspace{0.1cm}
	\subfloat[][\label{fig:cutting:Second}]{\includegraphics[width=0.49\textwidth]{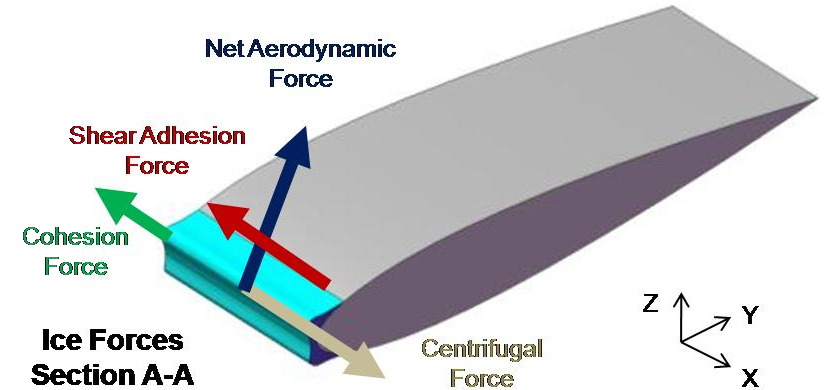}}
	\caption{Iterative cutting visualization \cite{brouwers2010experimental}: \textbf{(a)} initial partitioning, \textbf{(b)} ice piece analysed}
	\label{fig:cutting}
\end{figure}

The aim of this procedure is to analyse the ice shape by dividing it through $N$-cutting planes. As shown in Figure \ref{fig:cutting:First}, the rotor blade is divided into N ice "pieces" by cutting plane normal to the blade axis.

\begin{figure}[h!]
	\centering	
	\subfloat[][\label{CodeVisual:First}]{\includegraphics[width=0.4\textwidth]{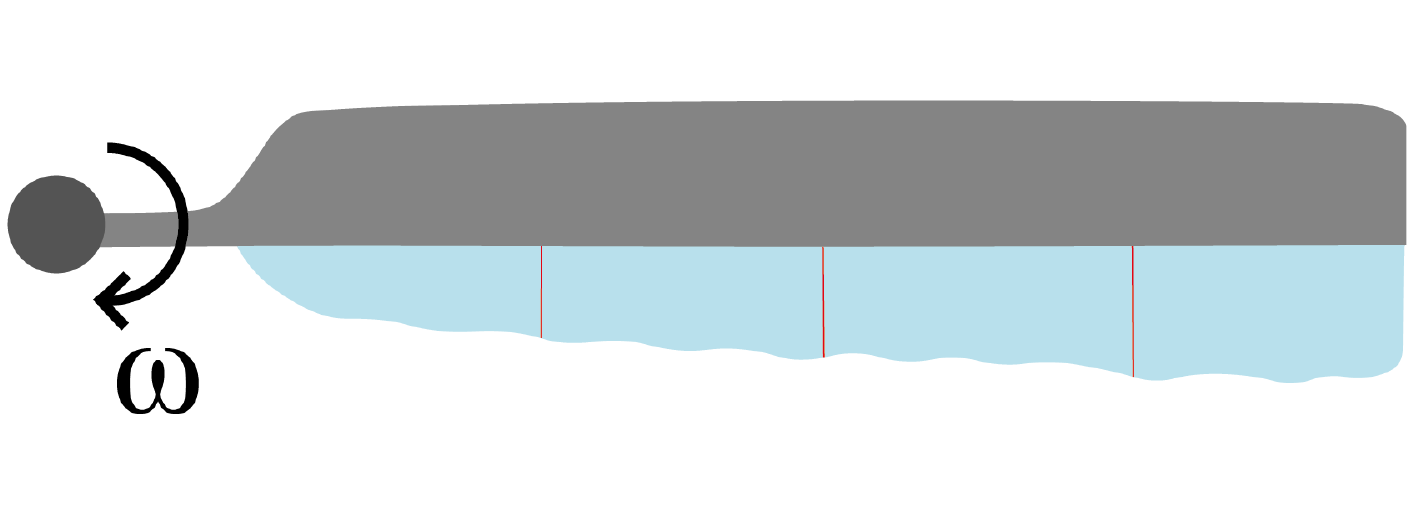}}
	\quad
	\subfloat[][\label{CodeVisual:Second}]{\includegraphics[width=0.4\textwidth]{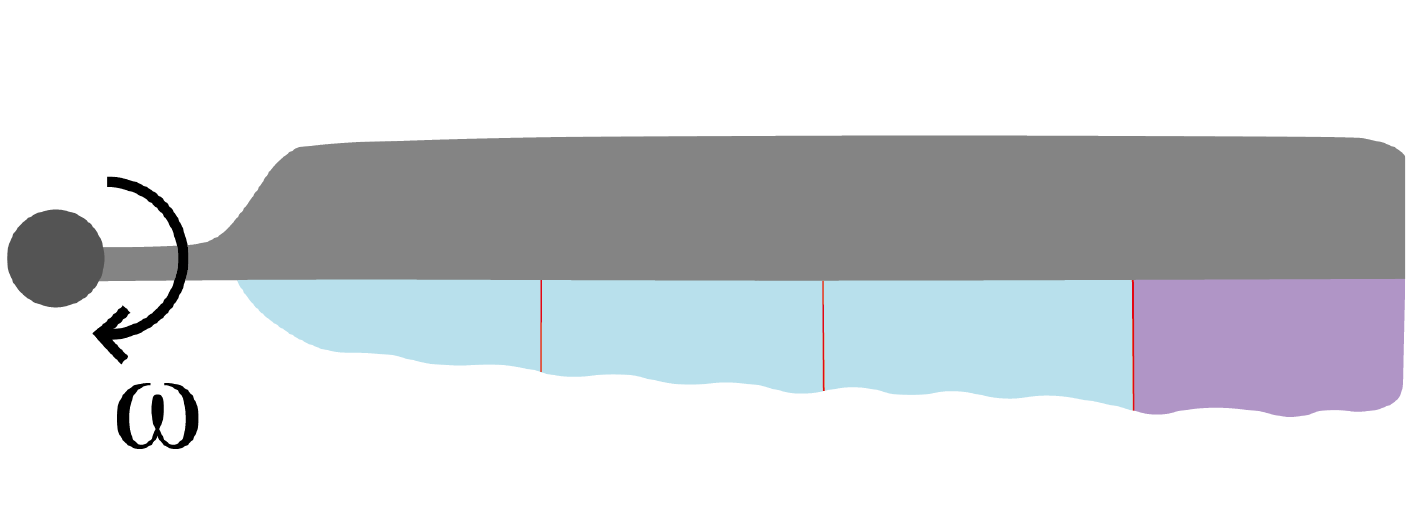}}
	\quad
	\subfloat[][\label{CodeVisual:Third}]{\includegraphics[width=0.4\textwidth]{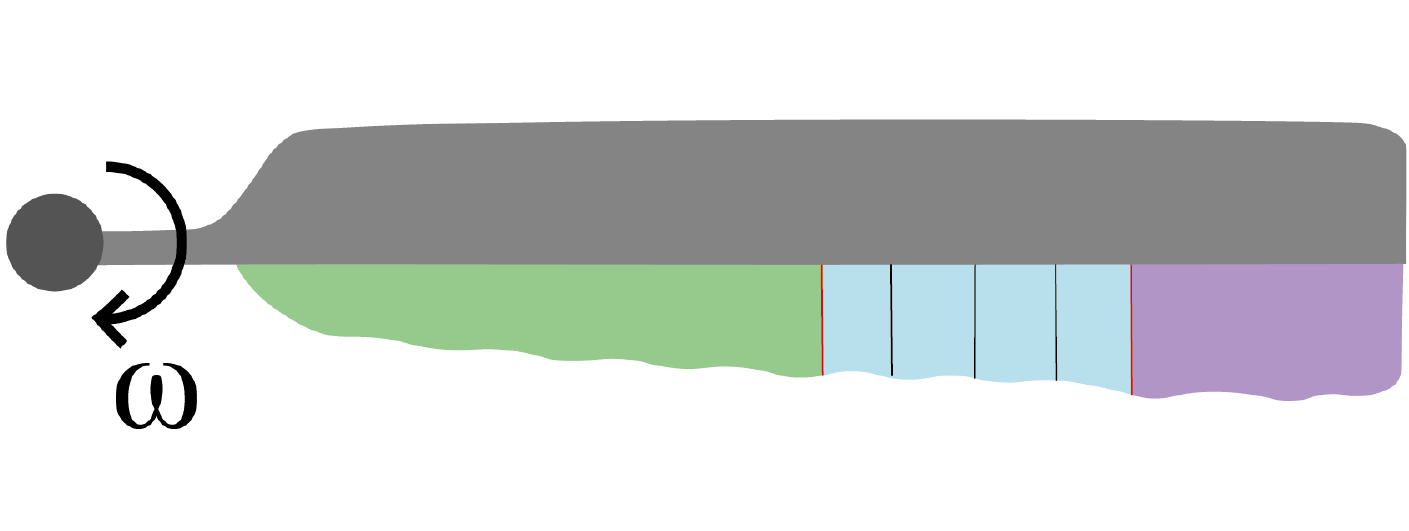}}
	\caption[Example of first 3 iterations]{Example of first 3 iterations: \textbf{(a)} initial partitioning, \textbf{(b)} not shedding case, \textbf{(c)} shedding case}
	\label{CodeVisual}
\end{figure}

Going from the tip to the root of the blade, the ice piece between two consecutive cutting planes is analysed by computing the centrifugal, cohesion and adhesion forces. The contribution of each ice piece is added to evaluate the shedding occurrence.
At each iteration, two cases can possibly occur:
\begin{enumerate}
	\item If shedding does not occur (Figure \ref{CodeVisual:Second}), the current ice portion is discarded (purple region), its contribution is considered for the following iterations and the next ice portion is analysed.
	\item If shedding occurs (Figure \ref{CodeVisual:Third}), the mesh of the ice that remains attached to the blade is discarded (green region), the current ice portion is set as the new computational domain (blue region), while keeping the contributions previously computed (purple region).
\end{enumerate} 

In this way it is possible to isolate the ice portion that will contain, if it exists, the shedding location, in order to find it with arbitrary precision through further adaptive domain-partitioning.

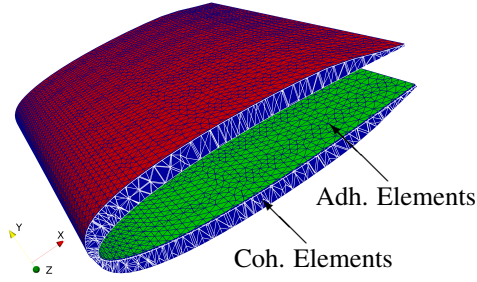
\begin{figure}[H]
	\centering
	\begin{overpic}[width=0.5\textwidth]{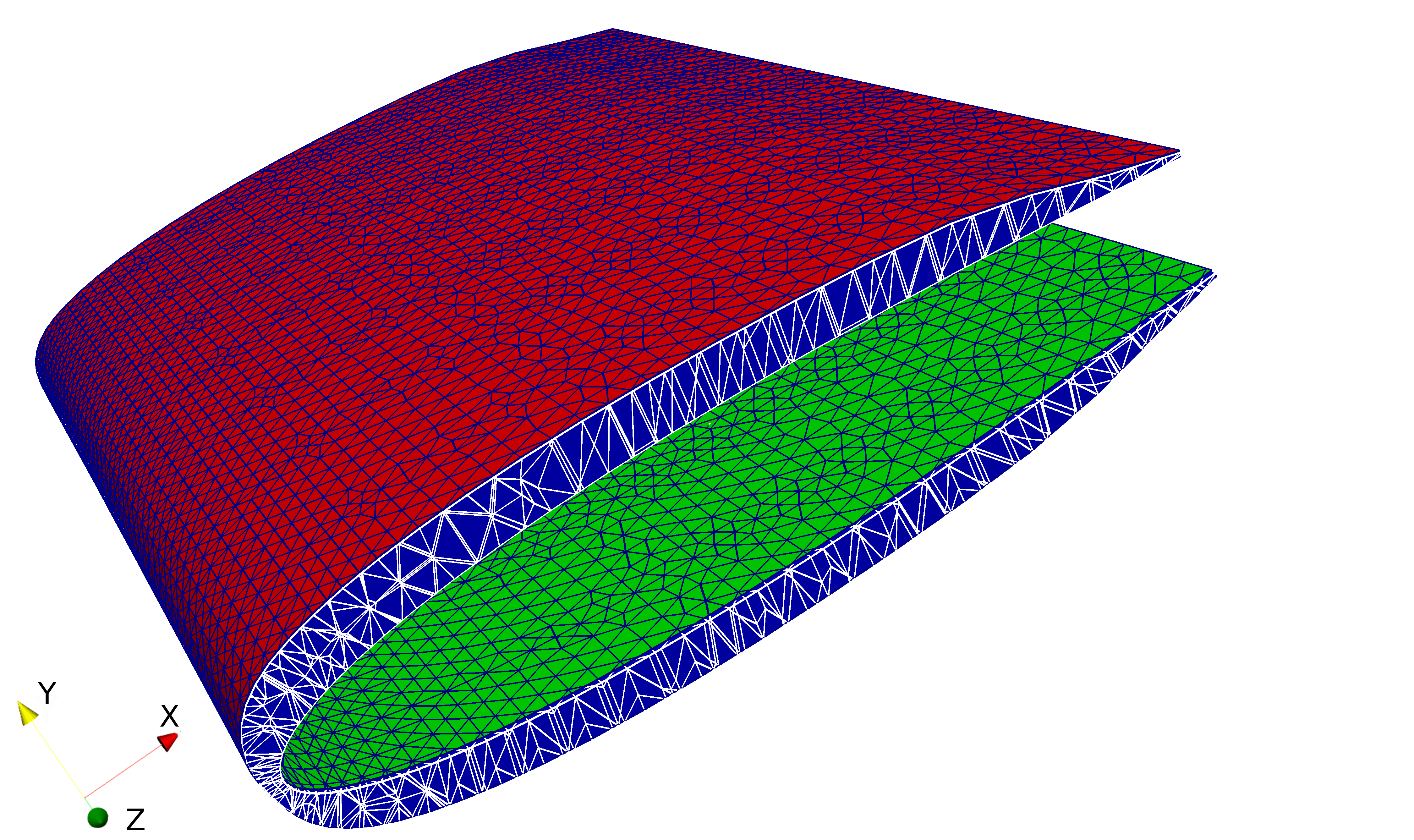}
		\put (48, 2) {\small Coh. Elements}
		\put (63,7) {\vector(-1, 1){9}}					
		\put (65, 15) {\small Adh. Elements}
		\put (80,20) {\vector(-1, 1){12}}
	\end{overpic}
	\caption{Subdivision of the elements of the mesh constituting the isolated ice piece}
	\label{fig:IceDivision}
\end{figure}

Forces are computed as sum of the contributions of the elements constituting the volume-mesh of the ice piece analysed. In particular, with reference to Figure \ref{fig:IceDivision}, the elements composing the ice mesh are gathered into four main categories:
\begin{itemize}
	\item \textit{Centrifugal elements}, which are the tetrahedrons constituting the mesh. They give their contribution ($F_{C,i}$) to the centrifugal force of the ice piece ($F_{C}$) 
	\begin{equation}\label{eq:CentrForce}
		F_{C} = \sum_{i=1}^{N_{e,vol}} F_{C,i} = \sum_{i=1}^{N_{e,vol}} m_i\Omega^2r_i
	\end{equation}

	where $N_{e,vol}$ is the number of centrifugal elements, $m_i$ and $r_i$ are respectively the mass of ice and the radial coordinate of the $i$-th centrifugal element, while $\Omega$ is the rotational velocity of the blade.
	\item \textit{Cohesion elements}, which are the volume elements that have intersections with the cutting plane (blue region). Their contribution ($F_{c,j}$) is summed up in the cohesion force of the ice piece ($F_{c}$)
	\begin{equation}\label{eq:CohForce}
		F_{c} = \sum_{j=1}^{N_{e,coh}} F_{c,j} = \sum_{j=1}^{N_{e,coh}} \sigma_{c}A_{c,j}
	\end{equation}

	where $N_{e,coh}$ is the number of cohesion elements, $A_{c,j}$ is the area of the intersection between the cutting plane and the $j$-th cohesion element and $\sigma_{c}$ is the cohesion strength of ice
	\item \textit{Adhesion elements}, which are the 3D surface elements that belong to the ice-blade interface (green region). They give their contribution ($F_{a,k}$) to the adhesion force of the ice piece ($F_{a}$)
	\begin{equation}\label{eq:AdhForce}
		F_{a} = \sum_{k=1}^{N_{e,adh}} F_{a,k} = \sum_{k=1}^{N_{e,adh}} \tau_{a}A_{a,k}
	\end{equation}

	where $N_{e,adh}$ is the number of adhesion elements, $A_{a,k}$ is the area of the $k$-th adhesion element and $\tau_a$ is the shear strength of ice
	\item \textit{Flow elements}, which are the 3D surface elements that belong to the ice-flow interface (red region). Their contribution is to the aerodynamic pressure of the ice piece, which, however, is not taken into account in the present work due to low values of Mach number involved and hence low values of the pressure coefficient \cite{scavuzzo1994influence}.
\end{itemize}

In previous works, Brouwers et al. \cite{brouwers2010experimental} used a surface mesh instead of a volume one. This forced them to analyse ice pieces with a span-wise constant area in order to compute the centrifugal force resulting in a use of a fine discretization. Using a volume-mesh, instead, and computing the centrifugal force as sum of the contributions of each element, it is possible to analyse the real shape of the ice as well as longer ice portions. This can lead to a speed-up of the convergence of a multi-step ice accretion analysis, in particular if the ice shedding does not occur.

\begin{figure}[H]
	\center
	\includegraphics[width=0.6\textwidth]{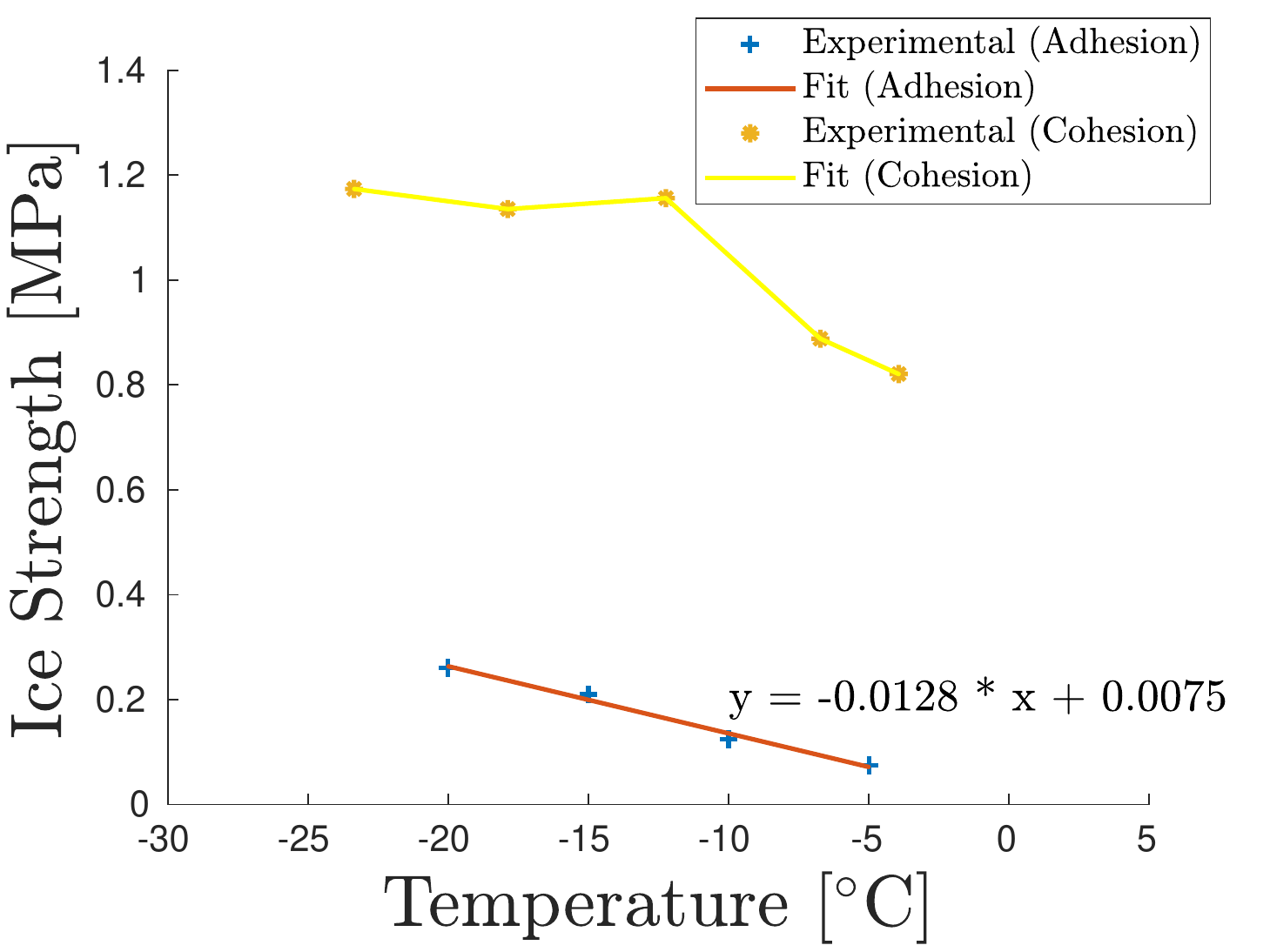}
	\caption{Mathematical formulations of the experimental models used in the present work for representing cohesion \cite{xian1989experimental} and adhesion \cite{IceAdhesionForShedding} strength of ice}
	\label{fig:ModelsUsed}
\end{figure}

From Equations (\ref{eq:CohForce}) and (\ref{eq:AdhForce}), cohesion and adhesion forces depend on the surface of the element considered and on two parameters, $\sigma_{c}$ and $\tau_{a}$. These parameters are, respectively, the cohesion and the adhesion strength of the ice and represent the values of the loads that have to be applied in order to cause a crack within the ice or its detachment from the supporting surface. In the open literature, the values found are very spread due to the different conditions in which test are performed and ice is generated. In particular, Loughborough \cite{brouwers2010experimental} used freezer ice, in which liquid water is allowed to slowly freeze to obtain the test material. This is not the situation for in-flight icing where supercooled droplets freezes upon impact with the blades (\textit{rime ice}) or flow over the surface to freeze downstream (\textit{glaze ice}). However, it is well-established the dependence of these parameters on the test air temperature, in particular they both decrease with an increase of the temperature. In Figure \ref{fig:ModelsUsed} are reported the values of the experimental results and the curves used for the fitting of these values for both adhesion and cohesion models taken into account in the present work.

\vspace{\baselineskip}

\noindent
In all simulations, the initial number of subdivisions is set to $N \ge 10$. However, in order to avoid too fine discretization, if no plane that satisfies the shedding conditions is found, a fitting procedure is performed on the three principal quantities analysed.

\subsection{Force-fitting procedure}\label{ch:Method:fitting}

During the iterative cutting it is possible that no cutting plane is found that satisfies the shedding conditions. Nevertheless ice shedding can still occur in between two consecutive cutting planes. Therefore, the aim of the \textit{Force-fitting procedure} is to avoid this situation, by searching for the shedding location along the whole blade.

\begin{figure}[h!]
	\centering	
	\subfloat[][\label{fig:SheddingCases:First}]{\includegraphics[width=0.45\textwidth]{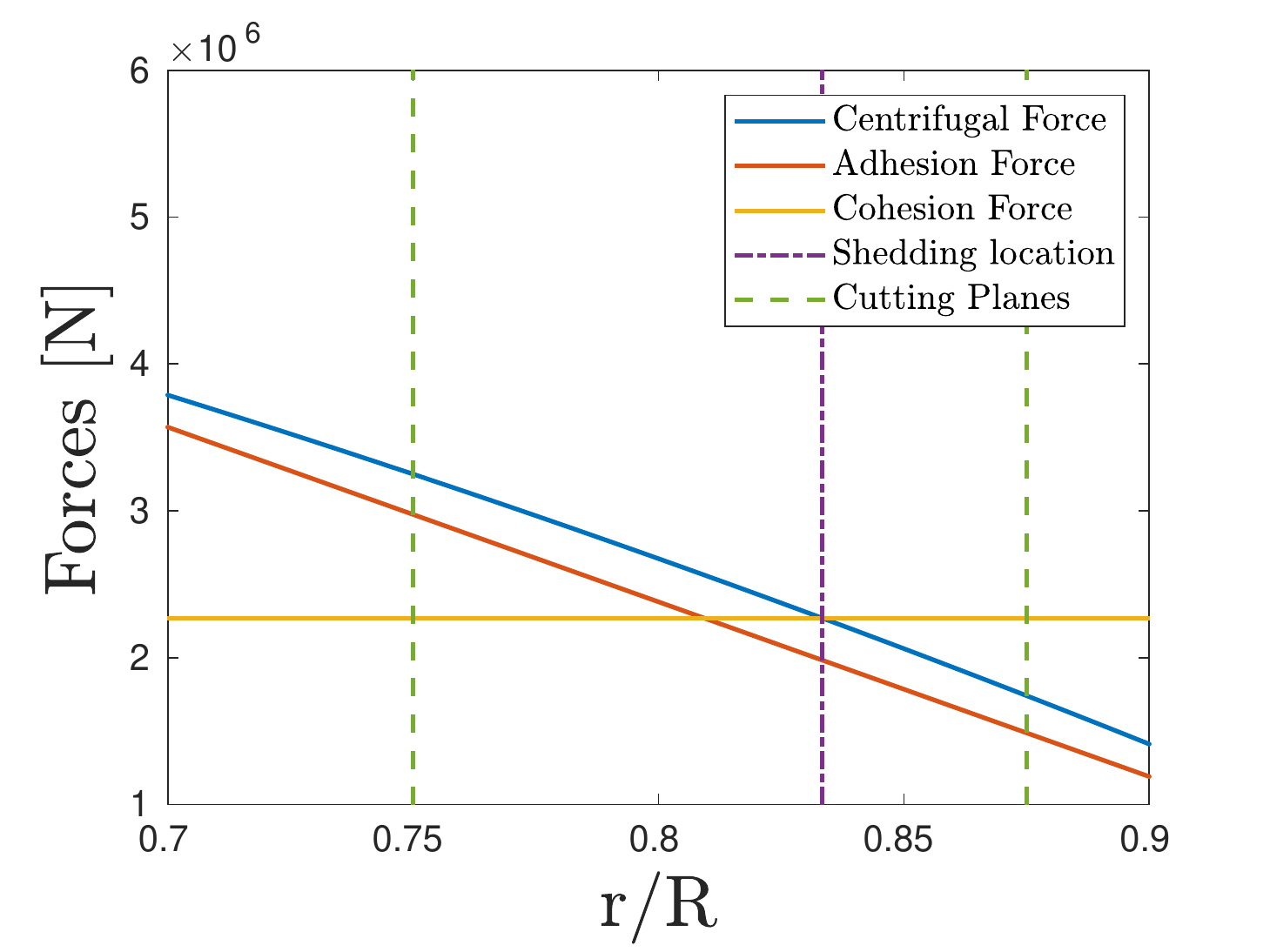}}
	\hspace{0.2cm}
	\subfloat[][\label{fig:SheddingCases:Second}]{\includegraphics[width=0.45\textwidth]{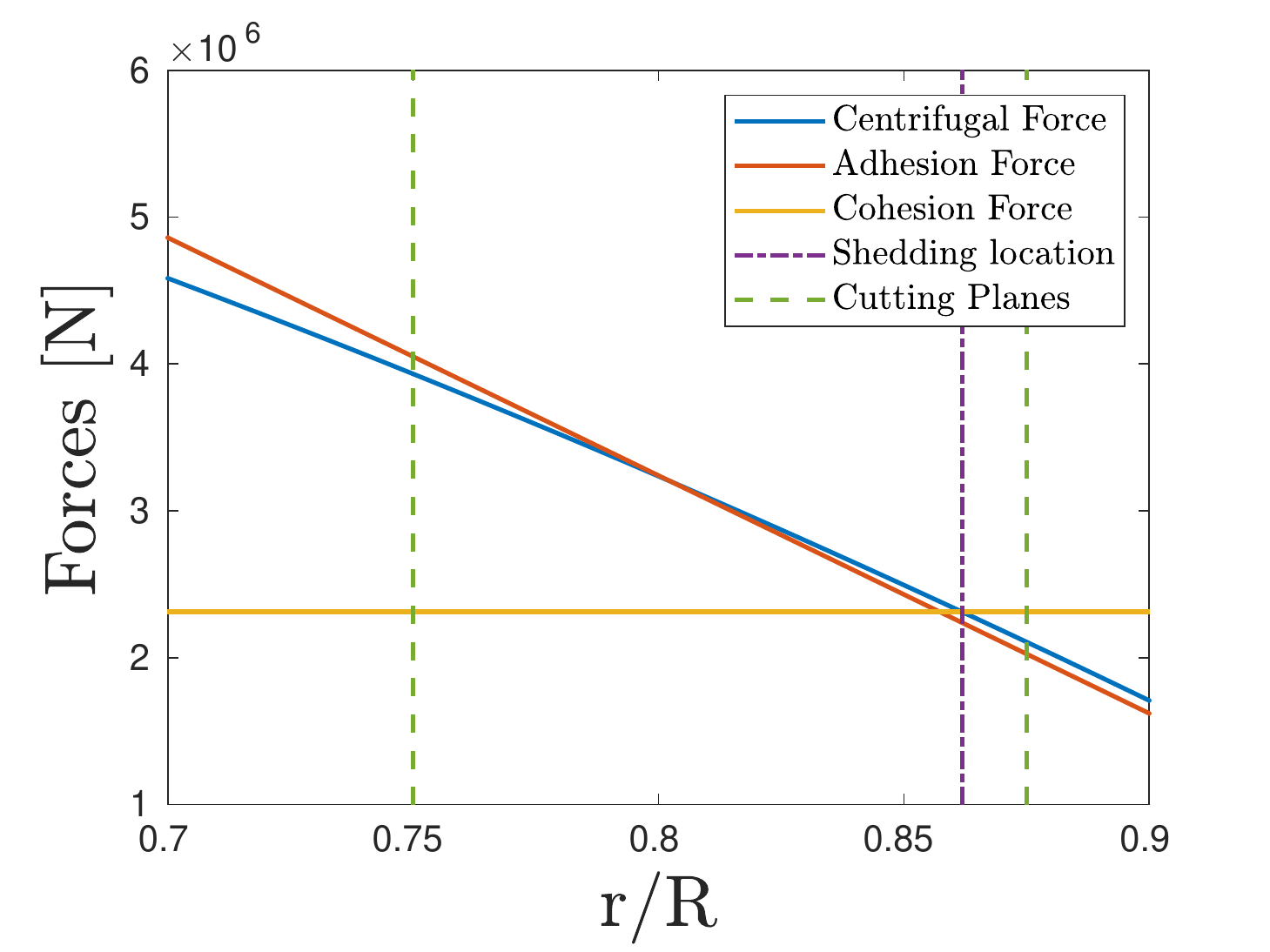}}
	\quad
	\subfloat[][\label{fig:SheddingCases:FirstDiff}]{\includegraphics[width=0.45\textwidth]{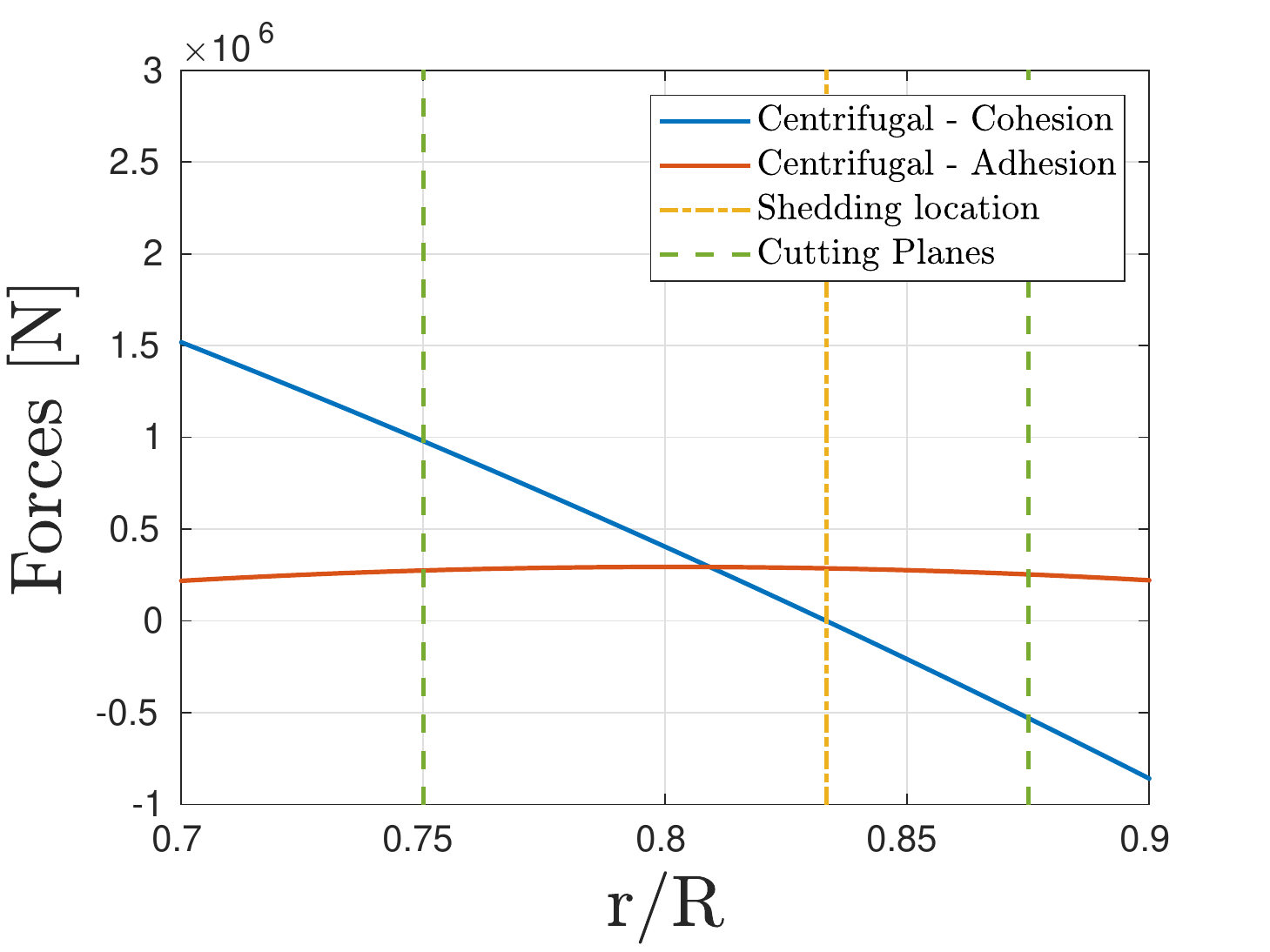}}
	\hspace{0.2cm}
	\subfloat[][\label{fig:SheddingCases:SecondDiff}]{\includegraphics[width=0.45\textwidth]{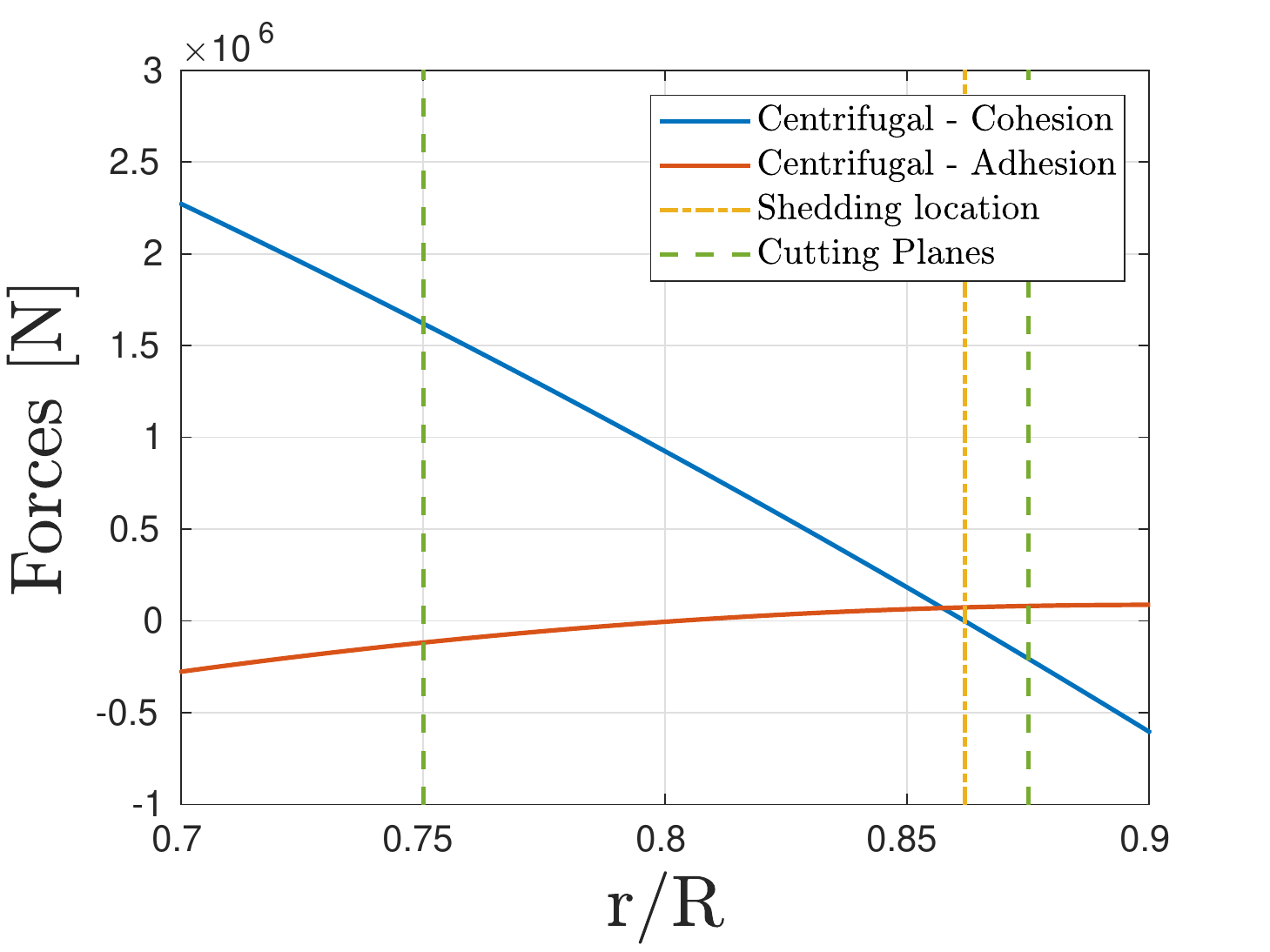}}
	\caption{Two possible cases of shedding occurrence}
	\label{fig:SheddingCases}
\end{figure}

\noindent
It can be reassumed with the following steps:
\begin{enumerate}
	\item Shedding analysis is performed. From tip to root the \textit{N}-cutting planes are investigated. 
	\item If shedding location is not found, a span-wise fitting is performed on centrifugal, adhesion and cohesion forces and the homonymous curves are obtained (Figure \ref{fig:SheddingCases}). 
	\item Using the fitted values and moving from tip to root, the radial coordinate $z_s$ at which shedding occurs is computed. If there is no such location the target point is now the one at which the curves of cohesion and centrifugal forces intersects. If again the point does not exist then the analysis is stopped and shedding does not occur. If it exists then this coordinate is $z_s$.
	\item The portion of ice that contains $z_s$ is isolated. The new domain to analyse is the one made out of the selected portion and its adjacent ones. The remaining volume towards the root of the blade is discarded, while keeping contributions of the portions placed on the tip side of the new domain. The new domain is divided in N portions and the shedding analysis restarts from step 1.
\end{enumerate}

The needs for the above fitting procedure can be better understood with reference to Figure \ref{fig:SheddingCases}. The curves reported represents the trends of centrifugal, adhesion and cohesion forces along the span of the blade for two different analytical cases, while the green dashed lines refer to the cutting planes and the purple dash-dotted lines refer to the shedding locations. Regarding the first case reported in Figure \ref{fig:SheddingCases:First}, there is a cutting plane at $r/R = 0.75$ that satisfies the shedding conditions. Thus the iterative cutting procedure is able to locate the shedding event through the procedure explained in Section \ref{ch:Method:Iterative}. On the other hand, the procedure developed without the force-fitting method would conclude that no shedding occurs for the case in Figure \ref{fig:SheddingCases:Second}. In fact, despite the presence of a shedding location at $r/R \approx 0.86$, there are no cutting planes that satisfy the shedding conditions. If, instead, the fitting procedure is taken into account, the coordinate $z_S$ will be found and, eventually, the shedding location will be determined.

\section{Results}

\begin{figure}[H]
	\center
	\includegraphics[width=0.7\textwidth]{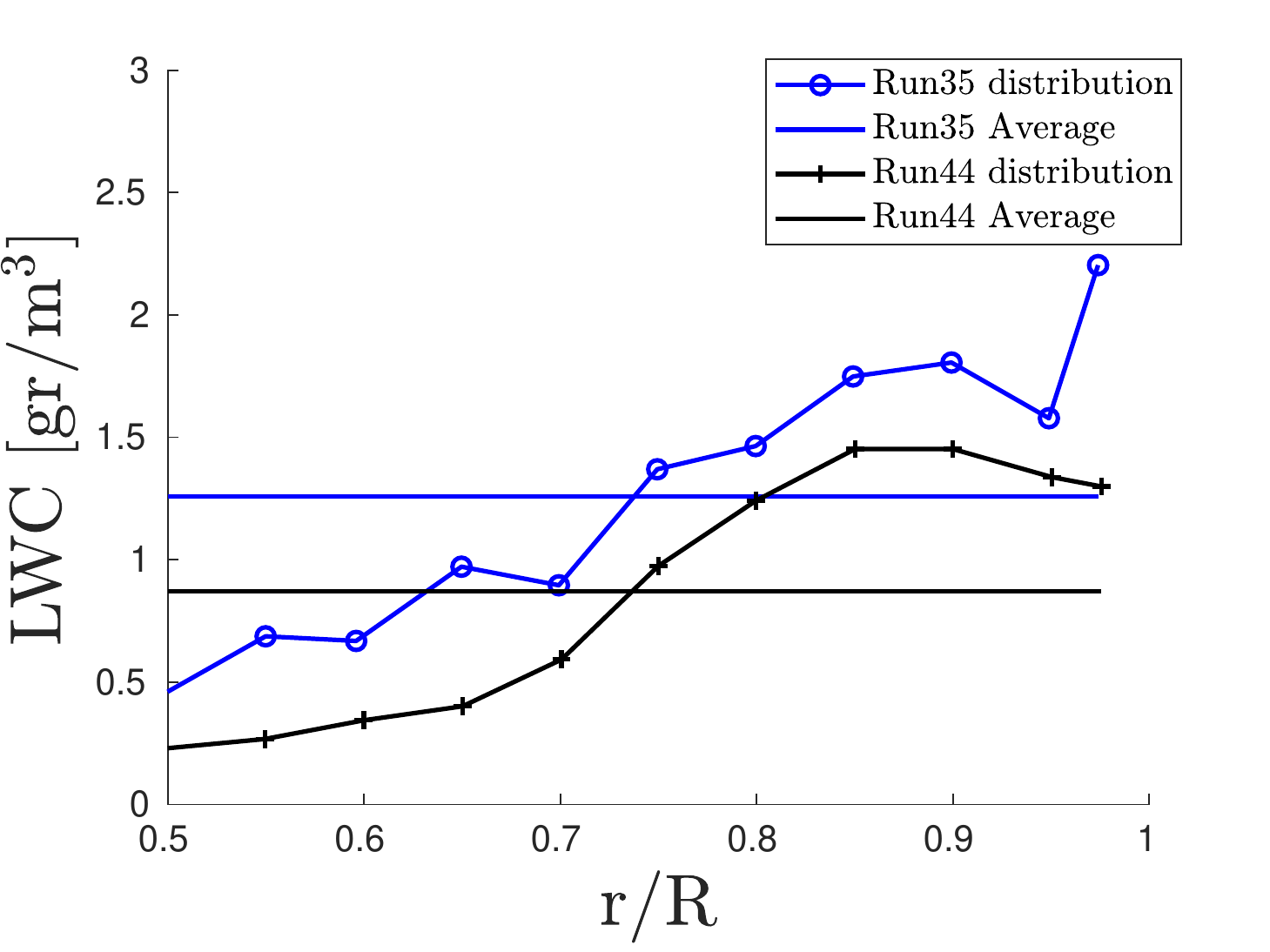}
	\caption{Span-wise distributions of LWC for ice accretion simulations \cite{brouwers2010experimentalThesis}}
	\label{fig:LWC}
\end{figure}

A Quasi-3D approach is used to obtain the ice shapes. 2D ice accretion simulations are performed at selected radial stations along the span of the blade. The 3D ice shape is then obtained by constant-area extrusion of the 2D geometries. 
The AERTS \cite{brouwers2010experimentalThesis} test cases are used to validate the ice shapes produced by the Quasi-3D approach as well as the ice shedding location and time. Due to absence of instruments to measure the Liquid Water Content (LWC) in the AERTS facility, the LWC has been computed by Brouwers et al. \cite{brouwers2010experimentalThesis} through empirical relations at radial locations along the blade. The values used for the test cases taken into account are reported in Figure \ref{fig:LWC}. The experimental runs are executed on a rotor with diameter of $2.36$\,m. The rotor tip speeds were limited to $600$\,RPM while the profile used for the blade sections was a NACA $0015$ with a constant chord span of $0.172$\,m. The blade rotor has a collective pitch of $2.5\,^{\circ}$ and a linear twist of $2.17\,^{\circ}$.

\begin{figure}[H]
	\centering	
	\subfloat[][\label{Run44ShapeComparison:06R}]{\includegraphics[width=0.46\textwidth]{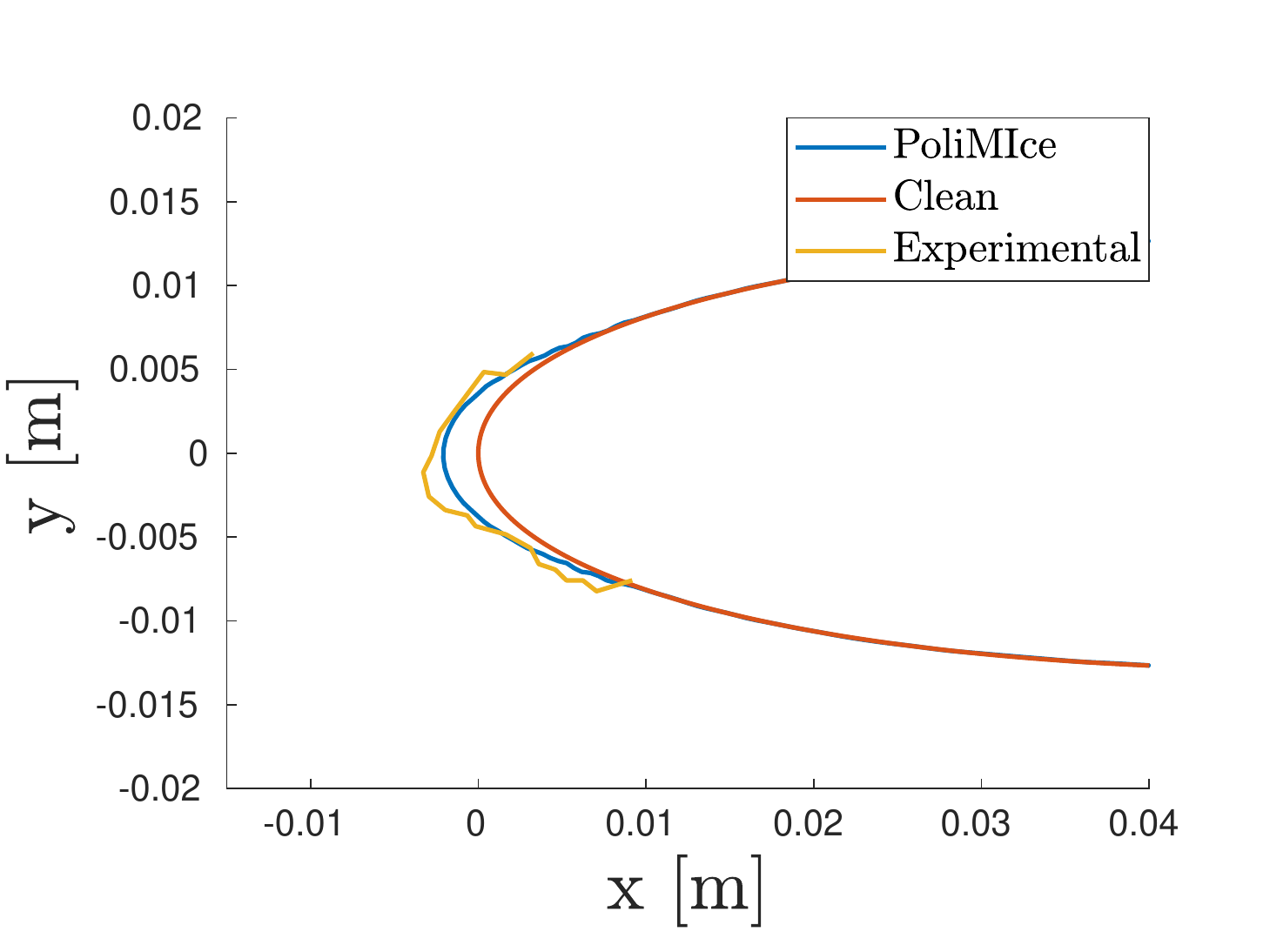}}
	\subfloat[][\label{Run44ShapeComparison:07R}]{\includegraphics[width=0.46\textwidth]{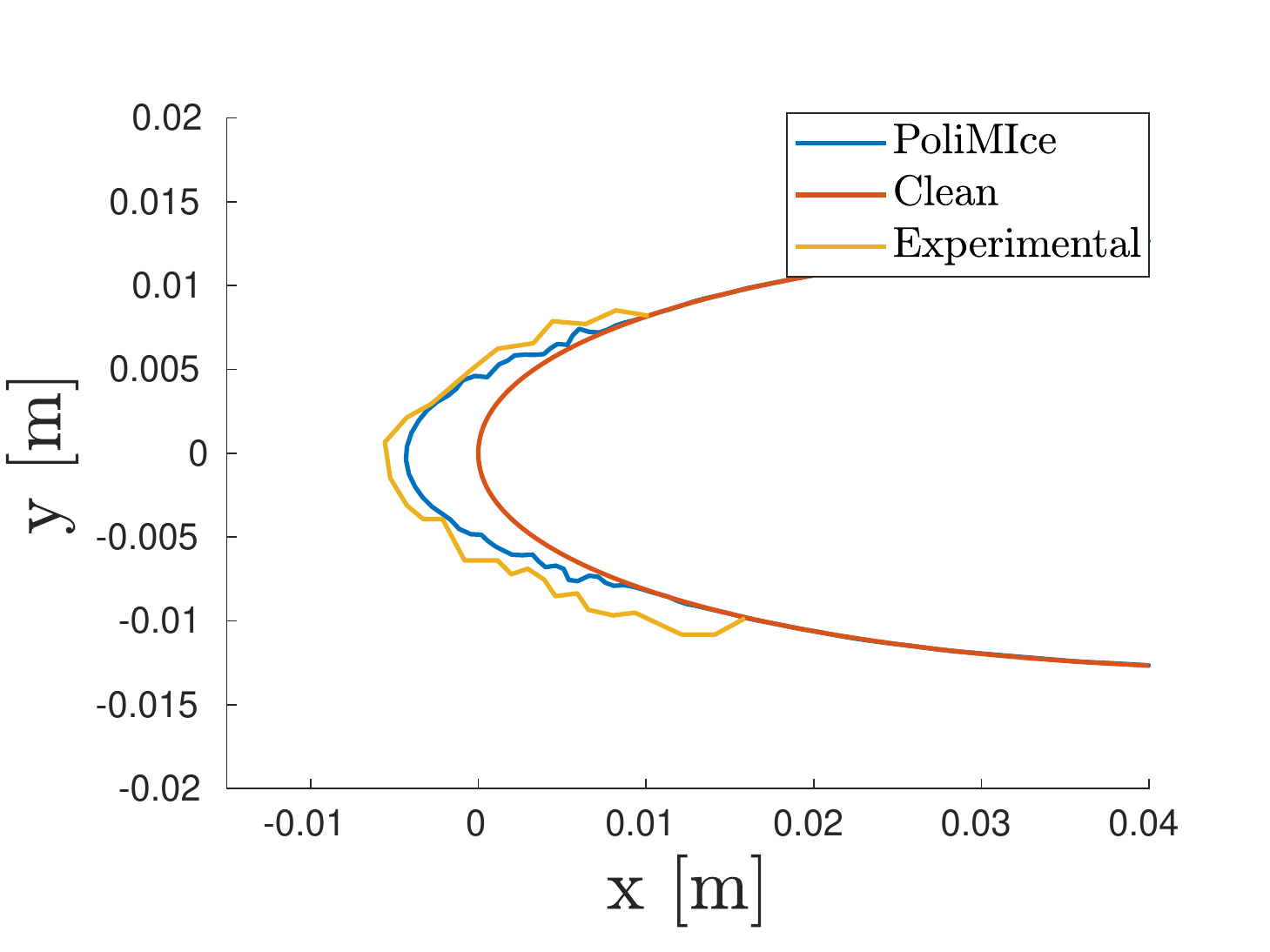}}
	\quad
	\subfloat[][\label{Run44ShapeComparison:08R}]{\includegraphics[width=0.46\textwidth]{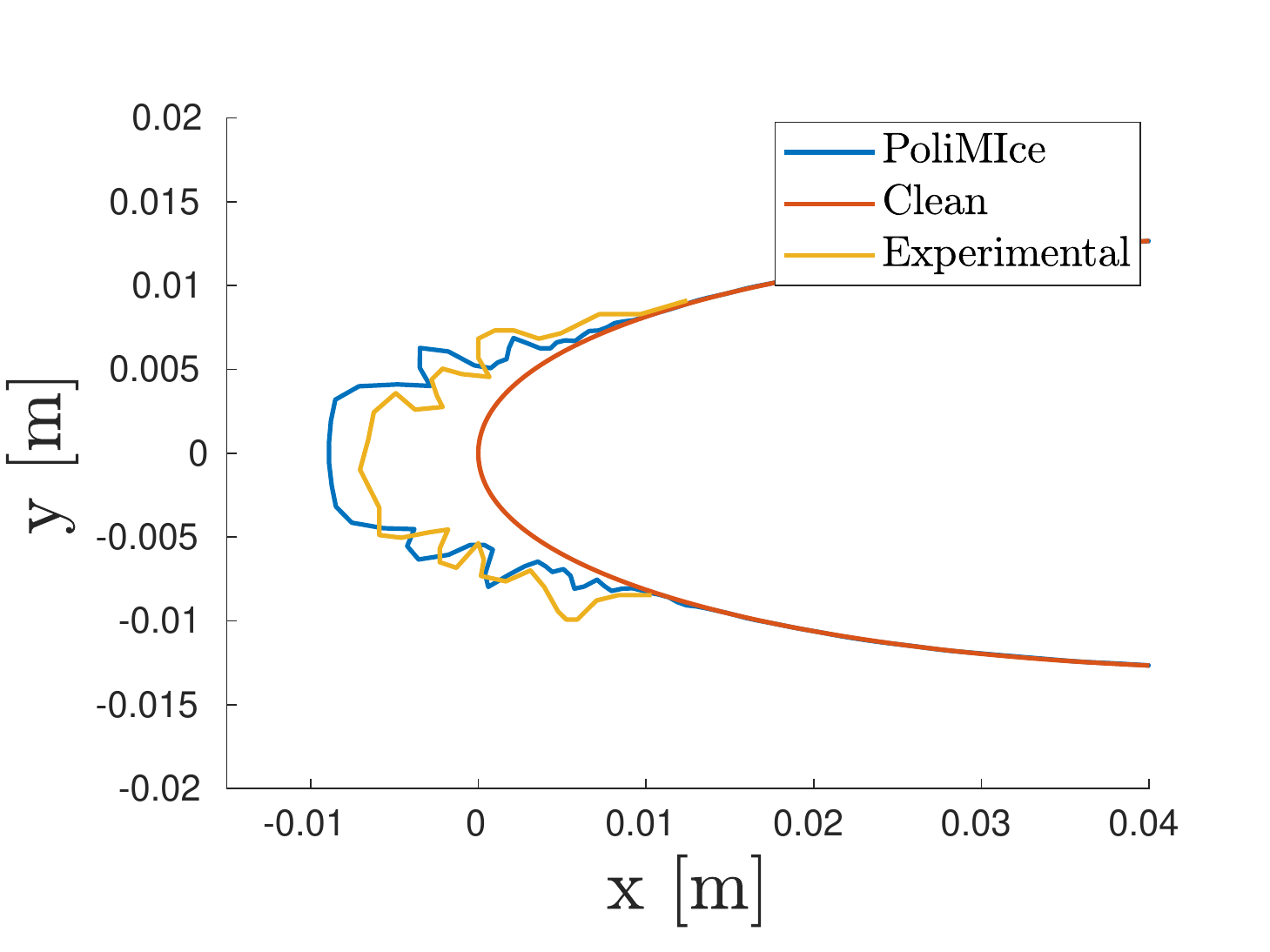}}
	\subfloat[][\label{Run44ShapeComparison:09R}]{\includegraphics[width=0.46\textwidth]{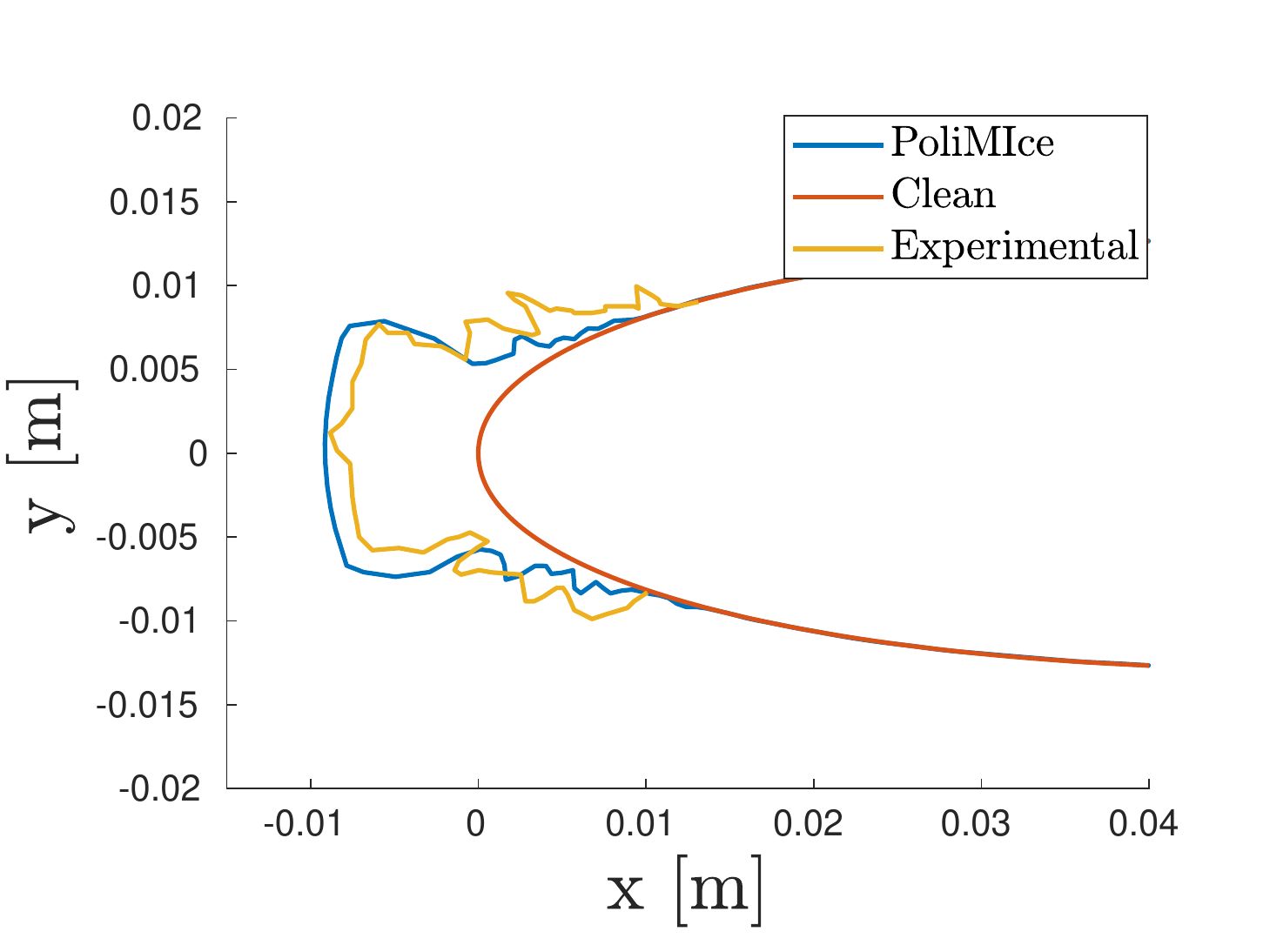}}
	\caption[Run 44 - Ice shape comparison at different stations along the blade]{Run 44 - Ice shape comparison at different stations along the blade: \textbf{(a)} $0.6$\,R, \textbf{(b)} $0.7$\,R, \textbf{(c)} $0.8$\,R, \textbf{(d)} $0.9$\,R}
	\label{fig:Run44ShapeComparison}
\end{figure}

\subsection{AERTS - Run 44}

Run 44 is used to assess the numerical sectional ice shapes against the experimental ones. The test conditions are reported in tab. \ref{tab:Run44}.

\begin{table}[h!]
	\begin{center}
		\begin{tabular}{l|r}
			\textbf{Parameter} & \textbf{Value}\\
			\hline
			Average temperature  & $-10.1$\,$^{\circ}$C \\
			Static Pressure & $101325$\,Pa \\
			MVD & $15$\,$\mu$m \\
			Ice accretion time & $180$\,s \\
			$\Delta\text{t}$ ice accretion & $30$\,s\\
			Ice shedding loc. & $\emptyset$
		\end{tabular}
		\caption{Test conditions AERTS\,--\,Run 44}
		\label{tab:Run44}
	\end{center}
\end{table}

Figure \ref{fig:Run44ShapeComparison} reports both the numerical and experimental ice shapes. The ice shapes obtained through Quasi-3D approach show good agreement with the experimental ones.

\subsection{AERTS - Run 35}

Run 35 is used for the validation of the proposed methodology for predicting ice shedding. The test conditions are reported in tab. \ref{tab:Run35}.

\begin{table}[h!]
	\begin{center}
		\begin{tabular}{l|r}
			\textbf{Parameter} & \textbf{Value}\\
			\hline
			Average temperature  & $-8.0$\,$^{\circ}$C \\
			Static Pressure & $101325$\,Pa \\
			MVD & $25$\,$\mu$m \\
			Ice accretion time & $360$\,s \\
			$\Delta\text{t}$ ice accretion & $40$\,s\\
			Ice shedding loc. & $0.83$\,R
		\end{tabular}
		\caption{Test conditions AERTS\,--\,Run 35}
		\label{tab:Run35}
	\end{center}
\end{table}

\begin{figure}[h!]
	\centering	
	\subfloat[][\label{fig:Run35Forces:160}]{\includegraphics[width=0.41\textwidth]{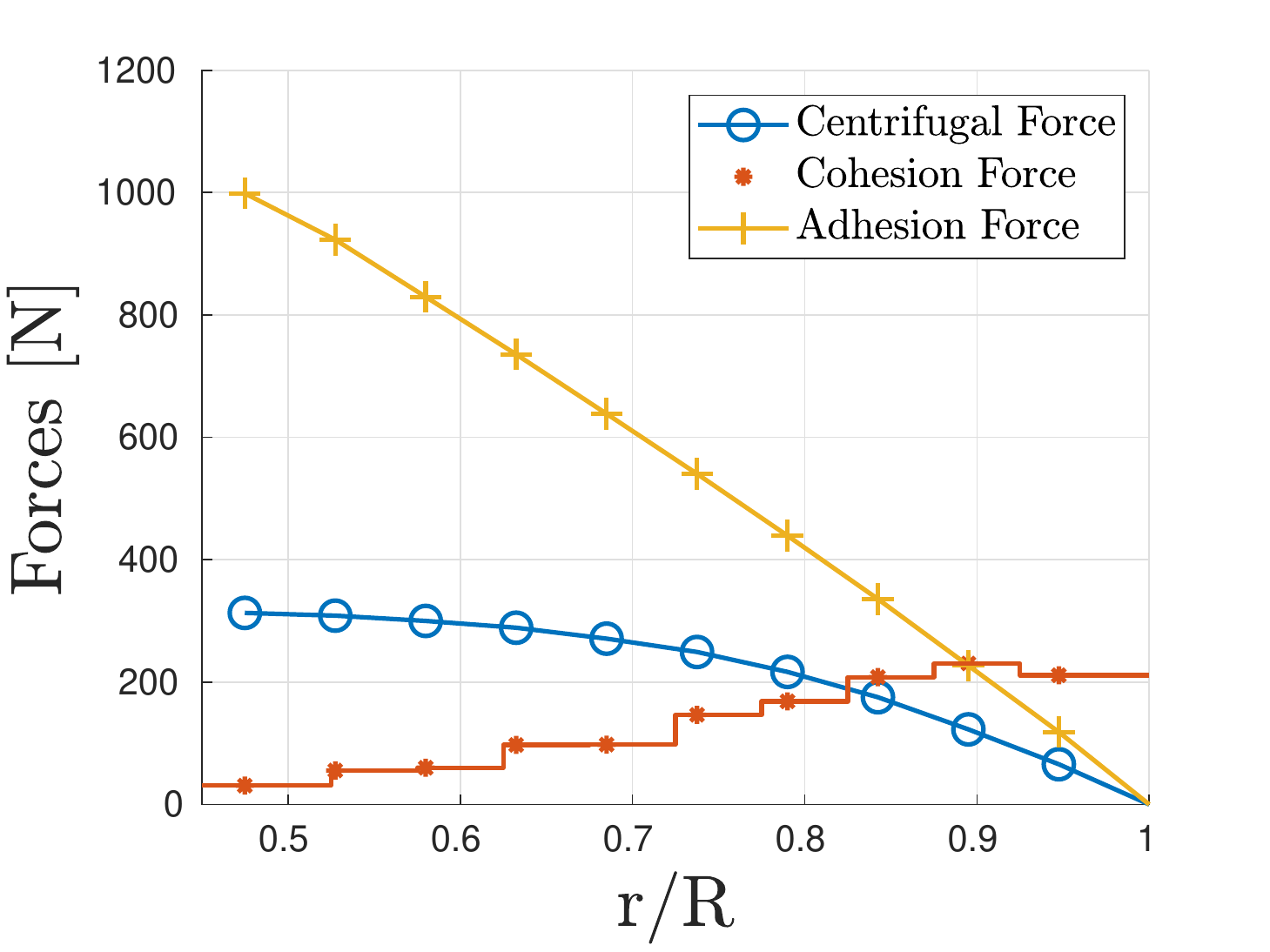}}
	\hspace{0.1cm}
	\subfloat[][\label{fig:Run35Forces:320}]{\includegraphics[width=0.41\textwidth]{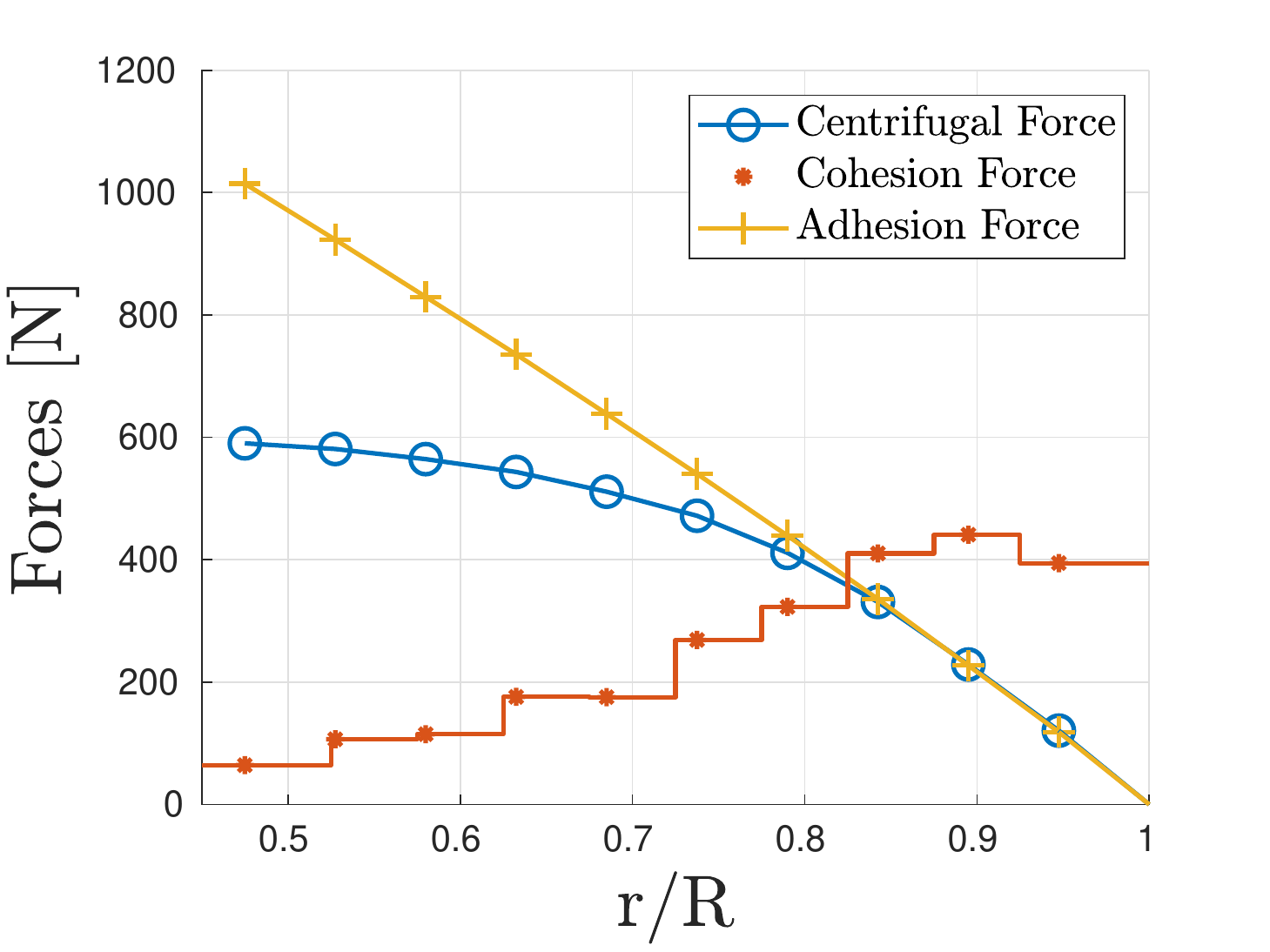}}
	\quad
	\subfloat[][\label{fig:Run35Forces:720}]{\includegraphics[width=0.43\textwidth]{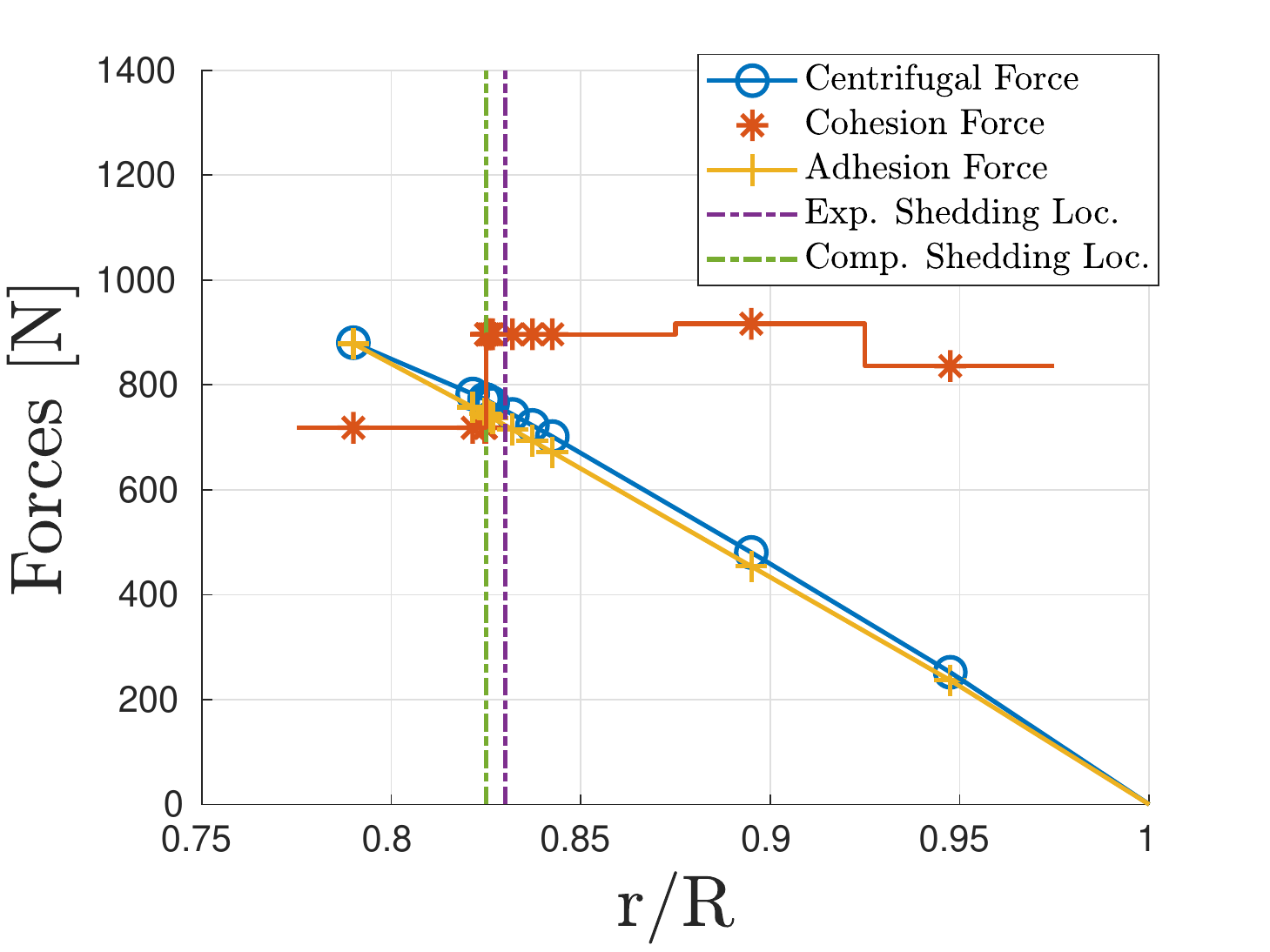}}
	\subfloat[][\label{fig:Run35Forces:720Zoom}]{\includegraphics[width=0.43\textwidth]{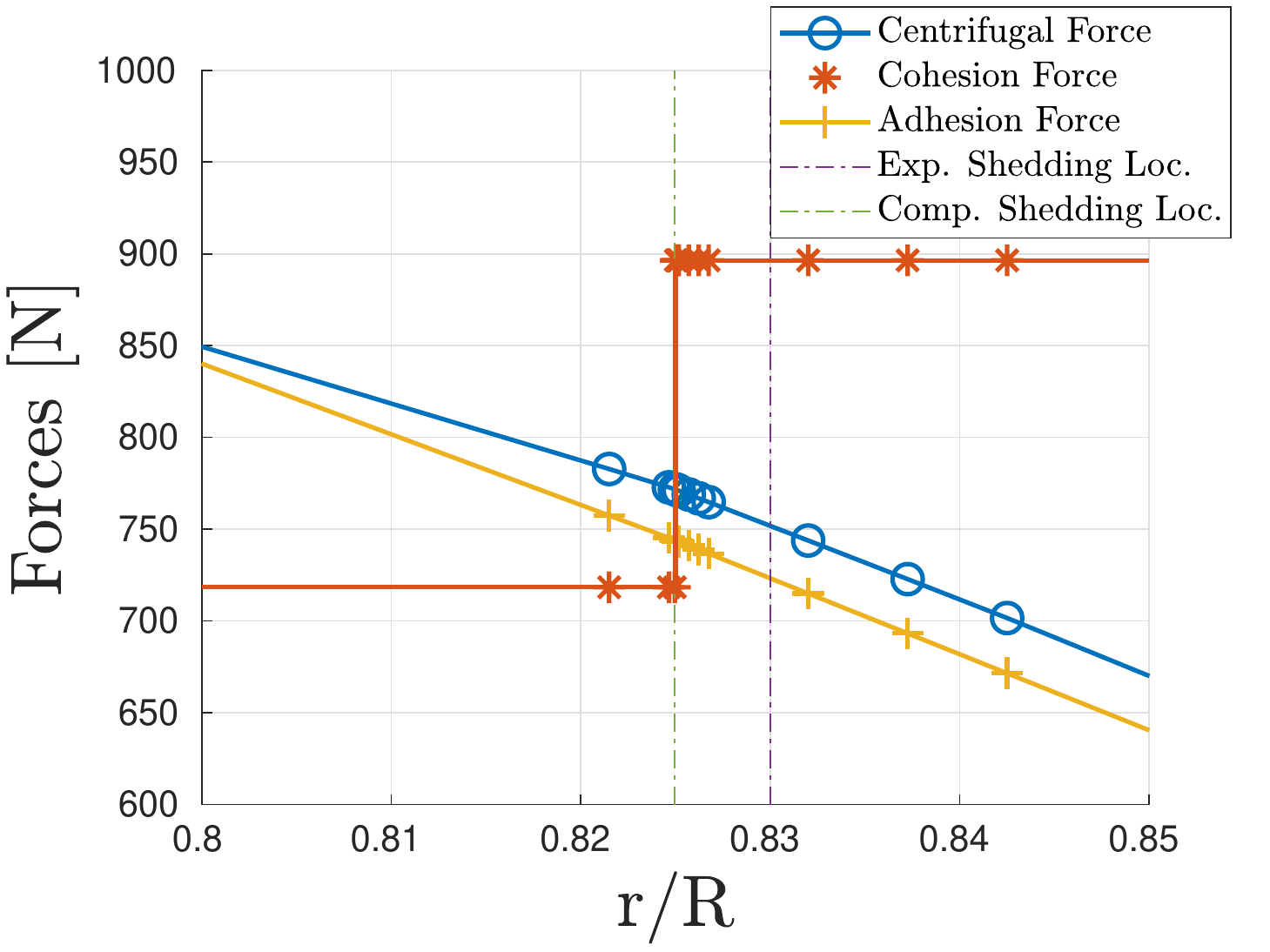}}
	\caption[Run 35 - Computed forces]{Run 35 - Computed forces: \textbf{(a)} at $160$\,s, \textbf{(b)} at $320$\,s, \textbf{(c)} at $720$\,s, \textbf{(d)} zoom at $720$\,s}
	\label{fig:Run35Forces}
\end{figure}

The span-wise trend of centrifugal, cohesion and adhesion forces is reported in Figure \ref{fig:Run35Forces} for three different time steps, respectively at $160$\,s, $320$\,s and $720$\,s. As it can be seen from these graphs, the adhesion force does not change its values during the accretion of the ice, whereas both cohesion and centrifugal forces increases due to the formation of new layers of ice. This, eventually, causes the shedding phenomenon, which, numerically, has occurred after $720$\,s and at a radial location of $0.825$\,R. Experimentally, ice shedding occurs after $320$\,s of ice accretion and at $0.83$\,R.

\newpage
\section{Conclusions}

A novel methodology for computing the shedding time and location for rotorcraft blades was presented. A Quasi-3D approach was used for computing the ice shapes along the blade. Good agreement has been achieved through this approach when comparing ice shapes obtained with experimental ones. An empirical evaluation of the shedding event is performed by analysing the centrifugal, cohesion and adhesion forces of increasingly larger ice pieces. By using a volume-mesh instead of a surface one it has been possible to analyse the real shape of the ice as well as longer ice pieces in order to make the present approach suitable for multi-step ice accretion simulations. The robustness of the method has been increased by a fitting procedure operated on the aforementioned forces. Numerical simulations are unable to compute correctly the timing of ice shedding. Further investigations on the models of ice adhesion and cohesion strength may increase the prediction accuracy of the proposed methodology.

\newpage
\section*{References}
\bibliographystyle{unsrt}
\bibliography{Bibliography}


\end{document}